\documentclass{amsart}
\usepackage{amsmath}
\usepackage{amsfonts}
\usepackage{enumerate}
\usepackage{amssymb}
\usepackage{mathbbol}
\usepackage{bm}

\DeclareMathOperator{\im}{Im} \DeclareMathOperator{\coker}{coker}
\DeclareMathOperator{\Ind}{Ind} \DeclareMathOperator{\Res}{Res}
\DeclareMathOperator{\co}{co}

\begin{document}
\theoremstyle{plain}
\newtheorem*{thm}{Theorem}
\newtheorem*{prop}{Proposition}
\newtheorem*{lem}{Lemma}
\newtheorem*{defn}{Definition}
\newtheorem*{ex}{Example}
\newtheorem*{cor}{Corollary}
\newtheorem*{quest}{Question}
\theoremstyle{remark}

\newcommand{\Md}[1]{\mathcal{M}(#1)}
\newcommand{\Pj}[1]{\mathcal{P}(#1)}
\newcommand{\G}[2]{\mathcal{G}_{#1}(#2)}

\title{On the Cartan map for crossed products and Hopf-Galois extensions}
\author{Konstantin Ardakov and Simon Wadsley}
\thanks{The first author thanks Christ's College, Cambridge for financial support}
\thanks{The second author was supported by EPSRC research grant EP/C527348/1}

\address{K. Ardakov: Department of Pure Mathematics, Hicks Building, Hounsfield Road, Sheffield S3 7RH, UK}
\email{K.Ardakov@shef.ac.uk}

\address{S. J. Wadsley: Department of Pure Mathematics and Mathematical Statistics, \newline Centre for Mathematical Sciences, Wilberforce Road, Cambridge CB3 0WB, UK}
\email{S.J.Wadsley@dpmms.cam.ac.uk} 

\begin{abstract}
We study certain aspects of the algebraic K-theory of Hopf-Galois
extensions. We show that the Cartan map from K-theory to G-theory of
such an extension is a rational isomorphism, provided the ring of
coinvariants is regular, the Hopf algebra is finite dimensional and
its Cartan map is injective in degree zero. This covers the case of
a crossed product of a regular ring with a finite group and has an
application to the study of Iwasawa modules.
\end{abstract}

\maketitle
\let\le=\leqslant  \let\leq=\leqslant
\let\ge=\geqslant  \let\geq=\geqslant

\section{Introduction}
\subsection{The Cartan map}\label{MainIntro} Recall that a ring is said to be \emph{right regular} if it is right Noetherian and every finitely generated right module has finite projective dimension. So any Noetherian ring of finite global dimension is necessarily regular.

One consequence of Quillen's celebrated Resolution Theorem is that the $K$-theory and the $\mathcal{G}$-theory of a right regular ring $B$ coincide \cite[Corollary 2 to Theorem 3]{Q}. More precisely, the Cartan map $K_i(B) \to \mathcal{G}_i(B)$ is an isomorphism for all $i\geq 0$.

Now if $G$ is a finite group and $A = B\ast G$ is a crossed product
then $A$ need not be regular, so the Resolution Theorem does not
apply. This is evident even in the simplest case when $B = k$ is a
field of characteristic $p > 0$, $p$ divides the order of $G$ and $A
= kG$ is the group algebra of $G$ --- in fact, in this case the
Cartan map $c : K_0(kG) \to \G{0}{kG}$ is not an isomorphism.

The purpose of this note is to prove the following
\begin{thm} Let $G$ be a finite group, let $k$ be a field, let $B$ be a right regular $k$-algebra and let $A = B \ast G$ be a crossed product. Then the Cartan map
\[c_i : K_i(A) \to \G{i}{A}\]
has torsion kernel and cokernel for all $i\geq 0$.
\end{thm}

Our original motivation for proving this theorem came from our study
of Euler characteristics of $p$-torsion Iwasawa modules \cite{AW}.
Using Theorem \ref{MainIntro} it is possible to extend the
definition of Euler characteristics to modules of infinite
projective dimension. See \cite[\S 8]{AW} for more details.

\subsection{Hopf-Galois extensions}
\label{HopfGalIntro} Given a Hopf algebra $H$ there is a notion of a
\emph{right $H$-Galois extension} of $k$-algebras $B \subset A$
\cite[8.1.1]{Mont}; the precise definition is given in (\ref{HG})
below. In the special case when the Hopf algebra $H$ is the group
algebra $kG$ of a group $G$, a result of Ulbrich \cite[Theorem
8.1.7]{Mont} states that an extension $B \subset A$ is right
$kG$-Galois if and only if $A$ is a strongly $G$-graded algebra
$A=\oplus_{g\in G}A_g$ with $B = A_1$. Thus any crossed product $A =
B \ast G$ a right $kG$-Galois extension.

Our main result (\ref{MainIntro}) is really a theorem about
Hopf-Galois extensions, as we make extensive use of the fact that
the comultiplication $\rho : A\to A\otimes H$ for a right
$H$-comodule algebra $A$ allows one to twist $A$-modules by
$H$-modules --- see (\ref{ModHG}). We actually prove the following
more general result.

\begin{thm}
Let $H$ be a finite dimensional Hopf algebra such that the Cartan
map $c : K_0(H) \to \G{0}{H}$ is injective. Let $B$ be a right
regular $k$-algebra and let $B \subset A$ be a right $H$-Galois
extension. There is a positive integer $m$ depending only on $H$
such that the kernel and cokernel of the Cartan map
\[c_i : K_i(A) \to \G{i}{A}\]
is killed by $m$ for all $i\geq 0$.
\end{thm}

Theorem \ref{MainIntro} then follows easily by appealing to the
classical theorem of Brauer \cite[Corollary 1 to Theorem 35]{Serre}
which states that for any field $k$ and any finite group $G$ the
Cartan map $c : K_0(kG) \to \G{0}{kG}$ is always injective. The
proof of Theorem \ref{HopfGalIntro} is given in (\ref{PfThm2}).

\subsection{Relaxing assumptions on $B$ and $H$}

It is natural to ask whether one can weaken the hypotheses in
Theorem \ref{HopfGalIntro}. The assumption that the Cartan map $c$
be injective is necessary: $H$ itself is always a right $H$-Galois
extension of the base field $k$ and Martin Lorenz has exhibited
examples of finite dimensional Hopf algebras $H$ for which $c$ is
not injective \cite[\S 4.2]{Lorenz}. As $K_0(H)$ is a torsion-free
abelian group it follows that the kernel of $c$ cannot be torsion in
this case.

Regarding the assumption on $B$, we ask the following

\begin{quest} Let $H$ be a finite dimensional Hopf algebra such that
the Cartan map $c : K_0(H) \to \G{0}{H}$ is injective and let $B
\subset A$ be a right $H$-Galois extension. Suppose that all Cartan
maps for $B$ have torsion kernel and cokernel. Does it follow that
all Cartan maps for $A$ also have torsion kernel and cokernel?
\end{quest}

\subsection{Conventions} Throughout $k$ will denote an arbitrary base field. All algebras, coalgebras, Hopf algebras, etc. are assumed to be defined over $k$. The unadorned tensor product $\otimes$ denotes $\otimes_k$. We will employ Sweedler's sumless notation for coalgebras and comodules: if $C$ is a coalgebra and $M$ is a right $C$-comodule with structure map $\rho:M\to M\otimes C$, then we will write
\[\rho(m) = m_0 \otimes m_1.\]
For an algebra $A$, $\Md{A}$ (respectively, $\Pj{A}$) will denote the category of all finitely generated (respectively, finitely generated projective) right $A$-modules. Unless specified otherwise, the term module will mean right module.

\subsection{Acknowledgements} The first author would like to thank
James Zhang, Ken Brown and Martin Lorenz for answering many of his
questions relating to this work during a conference in Shanghai.

\section{Hopf-Galois extensions}
\subsection{Comodule algebras} Let $A$ be an algebra and let $H$ be a Hopf algebra. The tensor product $A\otimes H$ is also an algebra with product given by the rule
\[(a\otimes b)\cdot (c\otimes d) = ac\otimes bd\qquad\mbox{ for all }a,c\in A, b,d \in H.\]
We say that $A$ is a \emph{right $H$--comodule algebra} if $A$ is a right $H$--comodule with structure map
\[\rho : A \to A\otimes H\]
such that $\rho$ is a map of algebras. Using Sweedler's sumless notation $\rho(a) = a_0\otimes a_1$ this second condition can be expressed as follows:
\[\rho(1) = 1\otimes 1 \qquad\mbox{ and }\qquad \rho(ab) = a_0b_0\otimes a_1b_1 \qquad\mbox{ for all   }a,b \in A.\]
Thus our definition agrees with the standard one given in \cite[4.1.2]{Mont}. We will usually write $B = A^{\co H}$ for the coinvariants of the coaction of $H$ on $A$:
\[B = A^{\co H} = \{a \in A : \rho(a) = a\otimes 1\}.\]
It is easy to verify that $B$ is always subalgebra of $A$, whenever $A$ is a right $H$--comodule algebra.

\subsection{Hopf-Galois extensions}\label{HG} Let $A$ be a right $H$--comodule algebra with coinvariants $B$. The extension $B \subset A$ is said to be a \emph{right $H$--Galois extension} \cite[8.1.1]{Mont} if the canonical map
\[\begin{array}{cccc}
\beta : & A\otimes_BA &\to& A\otimes H \\
 & x\otimes y &\mapsto& xy_0 \otimes y_1
\end{array}
\]
is bijective.

\begin{ex} The comultiplication $\Delta : H \to H\otimes H$ turns $H$ into a right $H$-comodule algebra. It can be checked that the coinvariants in this case are just $k$ and that $k \subset H$ is a right $H$-Galois extension.
\end{ex}

We will need the following result of Kreimer and Takeuchi \cite[Theorem 8.3.3]{Mont}:
\label{HopfGal}
\begin{thm} Let $H$ be a finite dimensional Hopf algebra and let $B\subset A$ be a right $H$--Galois extension. Then $A$ is a finitely generated projective right $B$-module.
\end{thm}

\subsection{Modules over Hopf-Galois extensions}
\label{ModHG}
Let $A$ be a right $H$--comodule algebra. Then whenever $M$ is an $A$-module and $V$ is an $H$--module, the tensor product $M\otimes V$ is naturally an $A\otimes H$--module. Using the algebra map $\rho : A \to A \otimes H$, $M\otimes V$ becomes an $A$-module, called the \emph{twist} of $M$ by $V$. The action of $A$ on $M\otimes V$ is explicitly given by the rule
\[(m\otimes v)\cdot a = ma_0 \otimes va_1 \qquad\mbox{ for all }m\in M, v\in V, a\in A.\]

Let $B$ be the coinvariants of $A$. If $M$ is an $A$-module and $N$ is a $B$-module, then we can form the \emph{restricted module} $\Res^A_B(M)$ and the \emph{induced module} $\Ind^A_B(N) := N \otimes_B A$. These constructions are related as follows.

\begin{prop} Let $B \subset A$ be a right $H$--Galois extension, let $M$ be an $A$-module and let $N$ be a $B$-module. Then
\begin{enumerate}[{(}a{)}]
\item $\Ind^A_B(\Res^A_B(M)) \cong M \otimes H$, and
\item $\Ind^A_B(N\otimes V) \cong \Ind^A_B(N) \otimes V$ provided the antipode of $H$ is invertible.
\end{enumerate}
\end{prop}
\begin{proof}(a) Note that $\beta : A\otimes_BA \to A\otimes H$ is a map of $A$--$A$-bimodules:
\[\begin{array}{rl}
\beta(a\cdot(x\otimes y)\cdot b) &= \beta(ax\otimes yb) = ax(yb)_0\otimes (yb)_1 = \\
 &= axy_0b_0\otimes y_1b_1 = a\cdot (xy_0\otimes y_1)\cdot b = a\cdot\beta(x\otimes y)\cdot b.
\end{array}
\]
As $B$ is a right $H$-Galois extension, $\beta$ is an isomorphism, so
\[1\otimes \beta : M\otimes_A(A\otimes_BA) \to M\otimes_A(A\otimes H)\]
is an isomorphism of $A$-modules. Part (a) follows.

(b) The action of $B$ on $N\otimes V$ is given by the rule $(n\otimes v)\cdot b = nb\otimes v$. Let $\sigma^{-1}$ be the inverse of the antipode $\sigma$ of $H$ and define maps
\[\begin{array}{lll}
\varphi &:&(N\otimes V)\otimes_BA \to (N\otimes_BA)\otimes V \qquad \mbox{ and } \\
\psi &:&(N\otimes_BA)\otimes V \to (N\otimes V)\otimes_BA
\end{array}\]
by the rules $\varphi(n\otimes v\otimes a) = n\otimes a_0\otimes va_1$ and $\psi(n\otimes a\otimes v) = n \otimes v\sigma^{-1}(a_1)\otimes a_0$. It is straightforward to verify the following asserions:
\begin{itemize}
\item $\varphi$ and $\psi$ are well-defined,
\item $\varphi$ is a map of $A$-modules,
\item $\varphi\circ\psi = 1_{(N\otimes_BA)\otimes V}$, and
\item $\psi\circ\varphi = 1_{(N\otimes V)\otimes_BA}$.
\end{itemize}
For the last two statements, use the fact that $\sigma^{-1}(h_2)h_1
= \epsilon(h) = h_2\sigma^{-1}(h_1)$ for all $h\in H$. Part (b)
follows.
\end{proof}

\subsection{The trivial module}
\label{triv}
The counit $\epsilon: H \to k$ gives $k$ the structure of a right $H$-module, called the \emph{trivial module} and denoted by $\mathbb{1}$.

\begin{lem} Let $A$ be a right $H$-comodule algebra. Then for any $A$-module $M$,
\[M\otimes \mathbb{1} \cong M\]
as $A$-modules. In particular, $V\otimes \mathbb{1} \cong V$ for any $H$-module $V$.
\end{lem}
\begin{proof} The first part follows from the counit axiom $a_0\epsilon(a_1) = a$ for the coaction $\rho : A \to A\otimes H$, and the second part follows as $H$ is itself a right $H$-comodule algebra as in Example \ref{HG}.
\end{proof}

\section{$K$-theory}
\subsection{Exact categories}
Let $\mathcal{B}$ be a full additive subcategory of an abelian category $\mathcal{A}$. We say that $\mathcal{B}$ is an \emph{exact category} if it is closed under extensions and kernels of surjections. That is, for any short exact sequence
\[0 \to A \to B \to C \to 0\]
in $\mathcal{A}$, the following conditions hold:
\begin{itemize}
\item if $A,C \in \mathcal{B}$ then $B\in\mathcal{B}$,
\item if $B,C \in \mathcal{B}$ then $A\in\mathcal{B}$.
\end{itemize}

The canonical example of an exact category is the category $\Pj{A} \subseteq \Md{A}$ for any right Noetherian algebra $A$.

\subsection{$K$-groups}
Quillen constructed a sequence of functors $K_0, K_1, \cdots$ which associate an abelian group $K_i(\mathcal{B})$ to every exact category $\mathcal{B}$. $K_0(\mathcal{B})$ is commonly known as the \emph{Grothendieck group} of $\mathcal{B}$. If $A$ is a right Noetherian algebra we will use the following shorthand notation:
\begin{itemize}
\item $\G{i}{A} := K_i(\Md{A})$, and
\item $K_i(A) := K_i(\Pj{A})$.
\end{itemize}

We will need the following well-known result about products in $K$-theory.
\label{ProdsKth}
\begin{lem} Let $\mathcal{A}, \mathcal{B}, \mathcal{C}$ be exact categories and let $\boxtimes : \mathcal{A} \times \mathcal{B} \to \mathcal{C}$ be an exact functor. Then for all $i,j\geq 0$ there exist bilinear maps
\[K_i(\mathcal{A}) \times K_j(\mathcal{B}) \to K_{i+j}(\mathcal{C})\]
which are natural in the variables $\mathcal{A}, \mathcal{B}, \mathcal{C}$.
\end{lem}

\subsection{The category $\mathcal{C}$}
\label{CatC}
From now on we will make the following assumptions:
\begin{itemize}
\item $H$ is a finite dimensional Hopf algebra,
\item $B \subset A$ is a right $H$--Galois extension,
\item $B$ is right Noetherian.
\end{itemize}

Theorem \ref{HopfGal} implies that $A$ is finitely generated as a $B$-module, so $A$ is also right Noetherian. Moreover, the restriction of a finitely generated $A$-module is finitely generated over $B$.

Let $\mathcal{C}$ denote the full subcategory of $\Md{A}$ consisting of modules $M$ whose restriction is a projective  $B$-module:
\[\mathcal{C} = \{M \in \Md{A} : \Res^A_B(M) \in \Pj{B}\}.\]
Since any $P\in\Pj{A}$ is a direct summand of $A^n$ for some $n\geq 0$, the restriction $\Res^A_B(P)$ is a finitely generated projective $B$-module by Theorem \ref{HopfGal}. Writing $\mathcal{P} := \Pj{A}$ and $\mathcal{M} := \Md{A}$ we have the following chain of inclusions:
\[\mathcal{P} \subseteq \mathcal{C} \subseteq \mathcal{M}.\]

\begin{lem} $\mathcal{C}$ is an exact category.
\end{lem}
\begin{proof} This is straightforward.
\end{proof}

\subsection{Tensor products of modules} The following result will be crucial to what follows. Note that as $H$ is assumed to be finite dimensional, $\Md{H}$ consists precisely of the finite dimensional $H$-modules.

\begin{lem}
\label{TensProd} Let $M\in \mathcal{M}$ and let $V, W \in \Md{H}$. Then
\begin{enumerate}[{(}a{)}]
\item $V \otimes W \in \Md{H}$,
\item $M\otimes V \in \mathcal{M}$,
\item if $M\in\mathcal{C}$ then $M\otimes V \in \mathcal{C}$,
\item if $M\in \mathcal{P}$ then $M\otimes V \in \mathcal{P}$,
\item if $M\in \mathcal{C}$ and $V\in \Pj{H}$ then $M\otimes V \in \mathcal{P}$.
\end{enumerate}
\end{lem}
\begin{proof} (a) By Example \ref{HopfGal}, $H$ is a right $H$-comodule algebra so the tensor product $V\otimes W$ is an $H$-module as in (\ref{ModHG}). This module is clearly finite dimensional.

(b) The restriction $\Res^A_B(M\otimes V)$ is isomorphic to $\Res^A_B(M)\otimes V$ with the $B$-action explained in the proof of Proposition \ref{ModHG}. The latter module is just a direct sum of $\dim V$ copies of $\Res^A_B(M)$ and is hence finitely generated over $B$. Hence $M\otimes V$ is finitely generated over $A$.

(c) This also follows from the isomorphism $\Res^A_B(M\otimes V)\cong \Res^A_B(M)\otimes V$.

(d) As tensor product commutes with direct sums, it is enough to show that $A\otimes V \in \mathcal{P}$. Because $H$ is finite dimensional, the antipode of $H$ is invertible by a result of Larson and Sweedler \cite[Theorem 2.1.3(2)]{Mont} so
\[A\otimes V = \Ind^A_B(B)\otimes V \cong \Ind^A_B(B\otimes V)\]
by Proposition \ref{ModHG}(b). Since $B\otimes V$ is isomorphic to $B^{\dim V}$, we see that $A\otimes V$ is isomorphic to $A^{\dim V}$ and is hence a finitely generated projective $A$-module.

(e) Again as tensor product commutes with direct sums it is enough to show that $M\otimes H$ is projective whenever $M\in \mathcal{C}$. But Proposition \ref{ModHG}(a) implies that $M\otimes H$ is isomorphic to $\Ind^A_B(\Res^A_B(M))$, which is a finitely generated projective $A$-module because $\Res^A_B(M)$ is a finitely generated projective $B$-module.
\end{proof}

\subsection{Module structures on $K$-groups}\label{ModStr}
\begin{prop} $\G{0}{H}$ is a ring. Moreover, $K_i(\mathcal{M})$, $K_i(\mathcal{C})$ and $K_i(\mathcal{P})$ are naturally right $\G{0}{H}$-modules for all $i\geq 0$.
\end{prop}
\begin{proof} By Lemma \ref{ProdsKth} and Proposition \ref{TensProd}(a) there is a well-defined associative product
$\G{0}{H} \times \G{0}{H} \to \G{0}{H}$. Using Lemma \ref{triv} we see that $\G{0}{H}$ is a ring with identity element $[\mathbb{1}]$.

A similar argument using Lemma \ref{triv}, Lemma \ref{ProdsKth} and parts (b),(c) and (d) of Proposition \ref{TensProd} shows that $K_i(\mathcal{M})$, $K_i(\mathcal{C})$ and $K_i(\mathcal{P})$ are all unital right $\G{0}{H}$-modules.
\end{proof}

\begin{cor} $K_0(H)$ is a right $\mathcal{G}_0(H)$-module.
\end{cor}
\begin{proof} This follows from Example \ref{HG}.
\end{proof}

\subsection{The Cartan maps}\label{Cartan} The inclusions $\Pj{H} \subseteq \Md{H}$ and $\mathcal{P} \subseteq \mathcal{C} \subseteq \mathcal{M}$ induce maps on the $K$-groups usually called \emph{Cartan maps}. We will use the following names for these:
\begin{itemize}
\item $c : K_0(H) \to \G{0}{H}$,
\item $c_i : K_i(\mathcal{P}) \to K_i(\mathcal{M})$,
\item $\lambda_i : K_i(\mathcal{P}) \to K_i(\mathcal{C})$, and
\item $\mu_i : K_i(\mathcal{C}) \to K_i(\mathcal{M})$.
\end{itemize}
Thus $c_i = \mu_i\circ \lambda_i$. By Lemma \ref{ModStr}, each group
appearing above is a $\G{0}{H}$-module. Moreover, all the maps
listed above are maps of right $\G{0}{H}$-modules.

\subsection{Resolutions}\label{Res}
\begin{prop} Suppose that the ring of coinvariants $B$ is right regular. Then every module $M\in\mathcal{M}$ has a finite resolution by objects in $\mathcal{C}$.
\end{prop}
\begin{proof} As $A$ is right Noetherian, we can find a projective resolution
\[\cdots \stackrel{d_n}{\longrightarrow} P_n \stackrel{d_{n-1}}{\longrightarrow} P_{n-1} \longrightarrow \cdots \stackrel{d_0}{\longrightarrow} P_0 \longrightarrow M \longrightarrow 0\]
consisting of finitely generated projective $A$-modules. We saw in (\ref{CatC}) that the restriction of each $P_n$ to $B$ is also finitely generated projective. Hence, writing $f_n = \Res^A_B(d_n)$ we see that
\[\cdots \stackrel{f_n}{\longrightarrow} \Res^A_B(P_n) \stackrel{f_{n-1}}{\longrightarrow} \cdots \stackrel{f_0}{\longrightarrow} \Res^A_B(P_0) \longrightarrow \Res^A_B(M) \longrightarrow 0 \]
is a projective resolution of $\Res^A_B(M)$. As $B$ is right regular, Schanuel's Lemma \cite[7.1.2]{MCR} implies that the $B$-module $\im(f_n) = \Res^A_B(\im(d_n))$ is projective for some $n\geq 0$. Hence $\im(d_n) \in \mathcal{C}$ and
\[0 \longrightarrow \im(d_n) \longrightarrow P_n \stackrel{d_{n-1}}{\longrightarrow} P_{n-1} \longrightarrow \cdots \stackrel{d_0}{\longrightarrow} P_0 \longrightarrow M \longrightarrow 0\]
is the required finite resolution of $M$ by objects in $\mathcal{C}$.
\end{proof}

\begin{cor} If $B$ is right regular then the map $\mu_i : K_i(\mathcal{C}) \to K_i(\mathcal{M})$ is an isomorphism for each $i\geq 0$.
\end{cor}
\begin{proof}
This follows from the Resolution Theorem \cite[Corollary 1 to
Theorem 3]{Q}.
\end{proof}

\subsection{Proof of Theorem \ref{HopfGalIntro}}
\label{PfThm2}

As $H$ is finite dimensional, $H$ is right Artinian. Hence $K_0(H)$
and $\G{0}{H}$ are finitely generated free abelian groups of the
same rank. As $c : K_0(H) \to \G{0}{H}$ is injective, $c$ has
torsion cokernel. It follows that we can find $P,Q \in
\mathcal{P}(H)$ such that $[P] - [Q] = m[\mathbb{1}]$ inside
$\G{0}{H}$ for some $m > 0$ depending only on $H$.

We have expressed $c_i$ as a product of two maps $c_i = \mu_i\circ\lambda_i$ in (\ref{Cartan}). By Corollary \ref{Res}, $\mu_i$ is an isomorphism so it is sufficient to show that the kernel and cokernel of $\lambda_i$ are killed by $m$.

Fix $i$ and write $\lambda = \lambda_i$ and $\mu = \mu_i$. For any $V\in\mathcal{M}(H)$, let $\tau_V$ denote the action of $[V] \in \G{0}{H}$ on $K_i(\mathcal{P})$; thus
\[\tau_V := K_i(M \mapsto M\otimes V) : K_i(\mathcal{P}) \to K_i(\mathcal{P}).\]
Now if $P\in\mathcal{P}(H)$, then $M\otimes P \in \mathcal{P}$ for all $M \in \mathcal{C}$ by Proposition \ref{TensProd}(e), so the map
\[\theta_P := K_i(M \mapsto M\otimes P) : K_i(\mathcal{C}) \to K_i(\mathcal{P})\]
is well-defined. Moreover, $\theta_P$ satisfies $\theta_P \circ \lambda = \tau_P$ for all $P\in\mathcal{P}(H)$.

Now $(\theta_P - \theta_Q)\circ\lambda  = \tau_P - \tau_Q =
m\tau_{\mathbb{1}} = m 1_{K_i(\mathcal{P})}$ as endomorphisms of
$K_i(\mathcal{P})$, because $K_i(\mathcal{P})$ is a unital
$\G{0}{H}$-module. Hence the kernel of $\lambda$ is killed by $m$.

Finally $\coker(\lambda)$ is a $\G{0}{H}$-module by (\ref{Cartan})
and it is annihilated by the image of the Cartan map $c : K_0(H) \to
\G{0}{H}$ by Proposition \ref{TensProd}(e). As $m[\mathbb{1}]$ lies
in the image of $c$, we see that $\coker(\lambda)$ is killed by $m$,
as required. \qed


\begin{thebibliography}{99}
\bibitem{AW} K. Ardakov, S. J. Wadsley, \emph{$K_0$ and the
dimension filtration for $p$-torsion Iwasawa modules}, submitted.
\bibitem{Lorenz} M. Lorenz, \emph{Representations of finite dimensional Hopf algebras}, J. Algebra \textbf{188} (1997), 476-505.
\bibitem{MCR} J.C. McConnell, J.C. Robson, \emph{Noncommutative Noetherian rings}, Revised Edition, AMS Graduate Studies in Mathematics, vol. \textbf{30} (2001).
\bibitem{Mont} S. Montgomery, \emph{Hopf algebras and their actions on rings}, CBMS Conference proceedings, AMS (1993).
\bibitem{Q} D. Quillen, \emph{Higher algebraic $K$-theory I.}, Lecture Notes in Mathematics, vol. \textbf{341}, Springer (1973), 85-147.
\bibitem{Serre} J.-P. Serre, \emph{Linear representations of finite groups}, Graduate Texts in Mathematics \textbf{42} Springer (1977).
\end{thebibliography}
\end{document}